\documentclass[11pt]{article}

\usepackage{geometry,graphicx,latexsym,amssymb,verbatim}
\usepackage{tikz}
\usetikzlibrary{calc}

\newtheorem{thm}{Theorem}

\newtheorem{ob}[thm]{Observation}

\newcommand{\smallqed}{{\tiny ($\Box$)}}

\newcommand{\qed}{$\Box$}

\def\vertex(#1){\put(#1){\circle*{2}}}
\def\vertexo(#1){\put(#1){\circle{2}}}
\def\vert(#1){\put(#1){\circle*{1.5}}}
\def\verto(#1){\put(#1){\circle{1.5}}}
\def\lab(#1)#2{\put(#1){\makebox(0,0)[c]{#2}}}
\def \nF {n_{_F}}
\def \mF {m_{_F}}

\newcommand{\TR}[1]{\mbox{$\tau(#1)$}}

\newcommand{\cC}{{\cal C}}

\def \nmr {\begin{enumerate}}
\def \enmr {\end{enumerate}}

\def \tmz {\begin{itemize}}
\def \etmz {\end{itemize}}

\newenvironment{unnumbered}[1]{\trivlist
\item [\hskip \labelsep {\bf
#1}]\ignorespaces\it}{\endtrivlist}

\newcommand{\hyperedgefour}[4]{
	\hyperedgebig #1#2#3;
	\hyperedgebig #2#3#4;
	\hyperedgebig #3#4#1;
	\hyperedgebig #4#1#2;
}

\newcommand{\hyperedgebig}[9]{
	\pgfmathsetmacro\Done{sqrt((#4-#1)^2+(#5-#2)^2)}
	\pgfmathsetmacro\angleone{(#2>#5)*(180+asin((#4-#1)/ \Done)-asin((#1-#4)/ \Done))+asin((#1-#4)/ \Done)+asin((#3-#6)/\Done)}
	\pgfmathsetmacro\angleone{\angleone-360*(\angleone>0)-360*(\angleone>360)}
	\pgfmathsetmacro\Dtwo{sqrt((#7-#4)^2+(#8-#5)^2)}
	\pgfmathsetmacro\angletwo{(#5>#8)*(180+asin((#7-#4)/ \Dtwo)-asin((#4-#7)/ \Dtwo))+asin((#4-#7)/ \Dtwo)+asin((#6-#9)/\Dtwo)}
	\pgfmathsetmacro\angletwo{\angletwo-360*(\angletwo>0)-360*(\angletwo>360)}
	\draw ([shift=(\angletwo:#6)] #4,#5)--([shift=(\angletwo:#9)]#7,#8);
	\draw (#4,#5)+(\angleone:#6) arc(\angleone:\angletwo+360*(\angletwo<\angleone):#6);
}

\newcommand{\hyperedgetwo}[6]{
	\pgfmathsetmacro\Done{sqrt((#4-#1)^2+(#5-#2)^2)}
	\pgfmathsetmacro\angleone{(#2>#5)*(180+asin((#4-#1)/ \Done)-asin((#1-#4)/ \Done))+asin((#1-#4)/ \Done)+asin((#3-#6)/\Done)}
	\pgfmathsetmacro\angleone{\angleone-360*(\angleone>0)-360*(\angleone>360)}
	\draw ([shift=(\angleone:#3)] #1,#2)--([shift=(\angleone:#6)]#4,#5);
	\pgfmathsetmacro\Dtwo{sqrt((#1-#4)^2+(#2-#5)^2)}
	\pgfmathsetmacro\angletwo{(#5>#2)*(180+asin((#1-#4)/ \Dtwo)-asin((#4-#1)/ \Dtwo))+asin((#4-#1)/ \Dtwo)+asin((#6-#3)/\Dtwo)}
	\pgfmathsetmacro\angletwo{\angletwo-360*(\angletwo>0)-360*(\angletwo>360)}
	\draw ([shift=(\angletwo:#6)] #4,#5)--([shift=(\angletwo:#3)]#1,#2);
	\draw (#1,#2)+(\angletwo:#3) arc(\angletwo:\angleone+360*(\angleone<\angletwo):#3);
	\draw (#4,#5)+(\angleone:#6) arc(\angleone:\angletwo+360*(\angletwo<\angleone):#6);
}
\begin{document}
\bibliographystyle{plain}

\title{A characterization of hypergraphs
that achieve \\ equality in the Chv\'{a}tal-McDiarmid Theorem}
\author{Michael A. Henning
 and Christian L\"{o}wenstein \\
\\
Department of Mathematics \\
University of Johannesburg \\
Auckland Park, 2006 South Africa\\
E-mail: mahenning@uj.ac.za; christian.loewenstein@uni-ulm.de}

\date{}
\maketitle

\begin{abstract}
For $k \ge 2$, let $H$ be a $k$-uniform hypergraph on $n$ vertices
and $m$ edges. The transversal number $\tau(H)$ of $H$ is the minimum
number of vertices that intersect every edge. Chv\'{a}tal and
McDiarmid [Combinatorica 12 (1992), 19--26] proved that $\tau(H)\le (
n + \left\lfloor \frac k2 \right\rfloor m )/ ( \left\lfloor
\frac{3k}2 \right\rfloor )$. When $k = 3$, the connected hypergraphs
that achieve equality in the Chv\'{a}tal-McDiarmid Theorem were
characterized by Henning and Yeo [J. Graph Theory 59 (2008),
326--348]. In this paper, we characterize the connected hypergraphs
that achieve equality in the Chv\'{a}tal-McDiarmid Theorem for $k =
2$ and for all $k \ge 4$.
\end{abstract}


\section{Introduction}

In this paper we continue the study of transversals in hypergraphs.
Hypergraphs are systems of sets which are conceived as natural
extensions of graphs.  A \emph{hypergraph} $H = (V,E)$ is a finite
set $V = V(H)$ of elements, called \emph{vertices}, together with a
finite multiset $E = E(H)$ of subsets of $V$, called
\emph{hyperedges} or simply \emph{edges}.

A $k$-\emph{edge} in $H$ is an edge of size~$k$. The hypergraph $H$
is said to be $k$-\emph{uniform} if every edge of $H$ is a $k$-edge.
Every (simple) graph is a $2$-uniform hypergraph. Thus graphs are
special hypergraphs.  The \emph{degree} of a vertex $v$ in $H$,
denoted by $d_H(v)$ or simply by $d(v)$ if $H$ is clear from the
context, is the number of edges of $H$ which contain $v$. The minimum
and maximum degrees among the vertices of $H$ is denoted by
$\delta(H)$ and $\Delta(H)$, respectively.

Two vertices $x$ and $y$ of $H$ are \emph{adjacent} if there is an
edge $e$ of $H$ such that $\{x,y\}\subseteq e$. The
\emph{neighborhood} of a vertex $v$ in $H$, denoted $N_H(v)$ or
simply $N(v)$ if $H$ is clear from the context, is the set of all
vertices different from $v$ that are adjacent to $v$.
Two vertices $x$ and $y$ of $H$ are \emph{connected} if there is a
sequence $x=v_0,v_1,v_2\ldots,v_k=y$ of vertices of $H$ in which
$v_{i-1}$ is adjacent to $v_i$ for $i=1,2,\ldots,k$. A
\emph{connected hypergraph} is a hypergraph in which every pair of
vertices are connected. A maximal connected subhypergraph of $H$ is a
\emph{component} of $H$. Thus, no edge in $H$ contains vertices from
different components.

If $H$ denotes a hypergraph and $X$ denotes a subset of vertices in
$H$, then $H-X$ will denote that hypergraph obtained from $H$ by
removing the vertices $X$ from $H$, removing all hyperedges that
intersect $X$ and removing all resulting isolated vertices, if any.
If $X = \{x\}$, we simply denote $H - X$ by $H - x$. We remark that
in the literature this is sometimes denoted by \emph{strongly
deleting} the vertices in $X$.

A subset $T$ of vertices in a hypergraph $H$ is a \emph{transversal}
(also called \emph{vertex cover} or \emph{hitting set} in many
papers) if $T$ has a nonempty intersection with every edge of $H$.
The \emph{transversal number} $\TR{H}$ of $H$ is the minimum size of
a transversal in $H$. A transversal of size~$\TR{H}$ is called a
$\TR{H}$-set.  Transversals in hypergraphs are well studied in the
literature (see, for example,~\cite{ChMc,HeLo12,HeYe08,HeYe10,HeYe13,LaCh90,ThYe07}). Chv\'{a}tal and McDiarmid~\cite{ChMc} established the following upper bound on
the transversal number of a uniform hypergraphs in terms of its order
and size.

\begin{unnumbered}{Chv\'{a}tal-McDiarmid Theorem.}
For $k \ge 2$, if $H$ is a $k$-uniform hypergraph on $n$ vertices
with $m$ edges, then
\[
\tau(H)\le \frac{n + \left\lfloor \frac k2\right\rfloor m}{\left\lfloor \frac{3k}2\right\rfloor}.
\]
\end{unnumbered}

As a special case of the Chv\'{a}tal-McDiarmid Theorem when $k = 3$,
we have that if $H$ is a $3$-uniform hypergraph on $n$ vertices with
$m$ edges, then $\tau(H) \le (n+m)/4$. This bound was independently
established by Tuza~\cite{Tu90} and a short proof of this result was
also given by Thomass\'{e} and Yeo~\cite{ThYe07}. The extremal
connected hypergraphs that achieve equality in the
Chv\'{a}tal-McDiarmid Theorem when $k = 3$ were characterized by
Henning and Yeo~\cite{HeYe08}. Their characterization showed that
there are three infinite families of extremal connected hypergraphs,
as well as two special hypergraphs, one of order~$7$ and the other of
order~$8$ .

Our aim in this paper is to characterize the connected hypergraphs
that achieve equality in the Chv\'{a}tal-McDiarmid Theorem for $k =
2$ and for all $k \ge 4$.
For this purpose we define two special families of hypergraphs.

\subsection{Special Families of Hypergraphs}
\label{S:special}

\begin{figure}[htb]
\begin{center}
\begin{tikzpicture}[scale=0.44]

\begin{scope}[xshift=-0.5cm,yshift=0cm]
\fill (0,0) circle (0.2cm); \fill (0,2) circle (0.2cm); \fill (0,4)
circle (0.2cm); \fill (0,6) circle (0.2cm);
\hyperedgetwo{0}{0}{0.7}{0}{6}{0.7}; \node at (0, -2) {\large $E_4$};
\end{scope}

\begin{scope}[xshift=8cm,yshift=0cm]
\fill (-2.7,0) circle (0.2cm); \fill (2.7,0) circle (0.2cm); \fill
(-3.45,1.3) circle (0.2cm); \fill (3.45,1.3) circle (0.2cm); \fill
(0,4.5) circle (0.2cm); \fill (0,6) circle (0.2cm); \hyperedgefour
{{-3.45}{1.3}{0.7}} {{-2.7}{0}{0.7}} {{2.7}{0}{0.6}}
{{3.45}{1.3}{0.6}}; \hyperedgefour {{0}{4.5}{0.7}} {{0}{6}{0.7}}
{{-3.45}{1.3}{0.6}} {{-2.7}{0}{0.6}}; \hyperedgefour {{2.7}{0}{0.7}}
{{3.45}{1.3}{0.7}} {{0}{6}{0.6}} {{0}{4.5}{0.6}}; \node at (0, -2)
{\large $T_4$};
\end{scope}

\begin{scope}[xshift=20cm,yshift=0cm]
\fill (-2.7,0) circle (0.2cm); \fill (2.7,0) circle (0.2cm); \fill
(-3.45,1.3) circle (0.2cm); \fill (1.95,1.3) circle (0.2cm); \fill
(3.45,1.3) circle (0.2cm); \fill (0,4.5) circle (0.2cm); \fill (0,6)
circle (0.2cm); \fill (-2.07,3.38) circle (0.2cm); \hyperedgefour
{{-3.45}{1.3}{0.7}} {{-2.7}{0}{0.7}} {{2.7}{0}{0.6}}
{{3.45}{1.3}{0.6}}; \hyperedgefour {{0}{4.5}{0.7}} {{0}{6}{0.7}}
{{-3.45}{1.3}{0.6}} {{-2.7}{0}{0.6}}; \hyperedgefour {{2.7}{0}{0.7}}
{{3.45}{1.3}{0.7}} {{0}{6}{0.6}} {{0}{4.5}{0.6}}; \node at (0, -2)
{\large $T_5$};
\end{scope}

\end{tikzpicture}
\end{center}
\vskip -0.6 cm \caption{The hypergraphs $E_4$, $T_4$, and $T_5$} \label{f:1}
\end{figure}
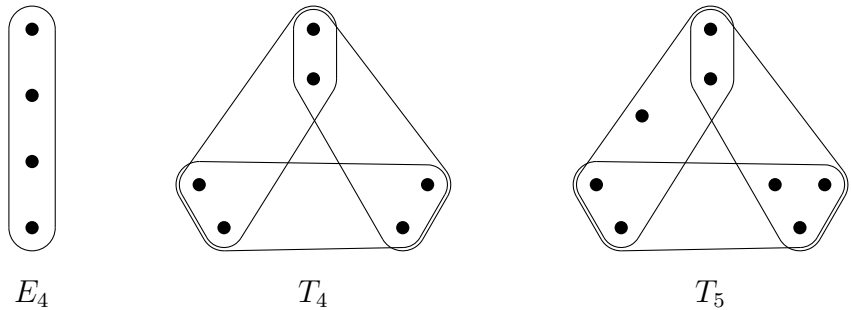

For $k \ge 2$, let $E_k$ denote the $k$-uniform hypergraph on $k$
vertices with exactly one edge. The hypergraph $E_4$ is illustrated
in Figure~\ref{f:1}.

For $k \ge 2$, a \emph{generalized triangle} $T_k$ is defined
as follows. Let $A$, $B$, $C$ and $D$ be vertex-disjoint sets of vertices with
$|A| = \left\lceil k/2 \right\rceil$, $|B| = |C| = \left\lfloor
k/2\right\rfloor$ and $|D| = \left\lceil k/2 \right\rceil -
\left\lfloor k/2 \right\rfloor$. In particular, if $k$ is even, the
set $D = \emptyset$, while if $k$ is odd, the set $D$ consist of a
singleton vertex. Let $T_k$ denote the $k$-uniform hypergraph with
$V(T_k) = A \cup B \cup C \cup D$ and with $E(T_k) =
\{e_1,e_2,e_3\}$, where $V(e_1) = A \cup B$, $V(e_2) = A \cup C$, and
$V(e_3) = B \cup C \cup D$. The hypergraphs $T_4$ and $T_5$ are
illustrated in Figure~\ref{f:1}.

\section{Main Result}

We shall prove:

\begin{thm}\label{t:main}
For $k = 2$ or $k \ge 4$, let $H$ be a connected $k$-uniform hypergraph on $n$ vertices and $m$ edges. Then,
\[
\tau(H) \le \frac{n + \left\lfloor \frac k2\right\rfloor m}{\left\lfloor \frac{3k}2\right\rfloor}
\]

\noindent with equality if and only if $H = E_k$ or $H = T_k$.
\end{thm}

We proceed as follows. We first recall some important results on edge
colorings of multigraphs in Section~\ref{S:color}. Thereafter we
establish a key theorem about matchings in multigraphs in
Section~\ref{S:match}. Finally in Section~\ref{S:proof} we present a
proof of Theorem~\ref{t:main} using an interplay between transversals
in hypergraphs and matchings in multigraphs.

\section{Edge Colorings of Multigraphs}
\label{S:color}

Let $G$ be a multigraph. An \emph{edge coloring} of $G$ is an
assignment of colors to the edges of $G$ such that adjacent edges
receive different colors. The minimum number of colors needed for an
edge coloring is called the chromatic index of the multigraph,
denoted $\chi'(G)$. The edge-multiplicity of an edge $e = uv$,
written $\mu(uv)$, is the number of edges joining $u$ and~$v$.
In his study of electrical networks, Shannon~\cite{Sh49} established
the following upper bound on the chromatic index of a multigraph.

\begin{unnumbered}{Shannon's Theorem.}
If $G$ is a multigraph, then $\displaystyle{ \chi'(G) \le
\left\lfloor 3\Delta(G)/ 2 \right\rfloor }$.
\end{unnumbered}

For $d \ge 2$, a \emph{Shannon multigraph of degree $d$} is a
multigraph on three vertices, with one pair of vertices joined by
$\left\lceil d/2 \right\rceil$ edges and the other two pairs joined
by $\left\lfloor d/2 \right\rfloor$ edges. Thus for fixed $d$, all
Shannon multigraphs of degree $d$ are isomorphic to the multigraph
$G$ with vertex set $V(G) = \{x,y,z\}$ and with $\mu(xy) =
\left\lfloor d/2 \right\rfloor$, $\mu(xz) = \left\lfloor d/2
\right\rfloor$ and $\mu(yz) = \left\lceil d/2 \right\rceil$. A
characterization of multigraphs achieving the upper bound in
Shannon's Theorem when the maximum degree is at least~$4$ was given
by Vizing~\cite{Vi65}.

\begin{unnumbered}{Vizing's Theorem.}
If $G$ is a connected multigraph with $\Delta(G)\ge 4$ and $\chi'(G)
= \left\lfloor 3\Delta(G)/2 \right\rfloor$, then $G$ contains a
Shannon multigraph of degree~$\Delta(G)$ as a submultigraph.
\end{unnumbered}

We remark the that maximum degree condition in Vizing's Theorem is
essential. For example, if $G$ is a connected multigraph with
$\Delta(G) = 2$ and $\chi'(G) = \left\lfloor 3\Delta(G)/2
\right\rfloor = 3$, then $G$ need not contain a Shannon multigraph of
degree~$\Delta(G)$ as a subgraph as may be seen by simply taking $G$
to be an odd cycle of length at least~$5$.

\section{Matchings in Multigraphs}
\label{S:match}

Let $G$ be a multigraph. Two edges in $G$ are \emph{independent} if
they are not adjacent in $G$. A set of pairwise independent edges of
$G$ is called a \emph{matching} in $G$, while a matching of maximum
cardinality is a \emph{maximum matching}. The number of edges in a
maximum matching of $G$ is called the \emph{matching number} of $G$
which we denote by $\alpha'(G)$. Our key matching theorem
characterizes connected multigraphs with small matching number
determined by Shannon's Theorem.

\begin{thm}\label{matching}
For $d \ge 4$, let $G$ be a connected multigraph of size $m$ with
$\Delta(G)\le d$. Then, $\alpha'(G) \ge m/\left\lfloor 3d/2
\right\rfloor$, with equality if and only if either $m = 0$ or $G$ is
a Shannon multigraph of degree~$d$.
\end{thm}
\textbf{Proof.} Let $\cC$ be an arbitrary edge coloring of the edges
of $G$ using $\chi'(G)$ colors. The matching number of $G$ is at
least the cardinality of a maximum edge color class in $\cC$, and so,
by Shannon's Theorem,
\[
\alpha'(G)\ge\frac m{\chi'(G)}\ge\frac m{\left\lfloor\frac{3\Delta(G)}2 \right\rfloor}
\ge\frac m{\left\lfloor\frac{3d}2 \right\rfloor},
\]

\noindent which establishes the desired lower bound. Suppose that
$\alpha'(G) = m/\left\lfloor 3d/2 \right\rfloor$ and $m \ge 1$. Then
we must have equality throughout the above inequality chain. Thus,
$\Delta(G)=d$, $\chi'(G)=\left\lfloor 3d/2 \right\rfloor$ and
$\alpha'(G) = m/\chi'(G)$. In particular, since $\cC$ is an arbitrary
$\chi'(G)$-edge coloring, the edge color classes in every
$\chi'(G)$-edge coloring have the same cardinality. Equivalently, the
edge color classes in $\cC$ are balanced. Since
$\chi'(G)=\left\lfloor 3d/2 \right\rfloor$, Vizing's Theorem implies
that $G$ contains a Shannon multigraph, $M$ say, of degree~$d$ as a
submultigraph.

If $d$ is even, then every vertex of $M$ has degree~$d$ in $M$. Since
$\Delta(G) = d$, the Shannon multigraph $M$ cannot be a proper
submultigraph of the connected multigraph $G$, implying that $G = M$.
Hence if $d$ is even, then $G$ is a Shannon multigraph of degree $d$.
Therefore we may assume that $d$ is odd, for otherwise the desired
result holds.

Since $d$ is odd, $d \ge 5$ and one pair of vertices in $M$ is joined
by $(d+1)/2$ edges and the other two pairs are joined by $(d-1)/2$
edges. Thus two vertices in $M$ have degree~$d$ in $M$ and one
vertex, $x$ say, of $M$ has degree~$d-1$ in $M$. Assume that $M$ is a
proper submultigraph of $G$. Since $\Delta(G) = d$, the vertex $x$ is
adjacent in $G$ to exactly one vertex $v \notin V(M)$. Since $xv$ is
a bridge in $G$, the edge $xv$ cannot belong to a submultigraph of
$G$ that is isomorphic to a Shannon multigraph of degree~$d$. Thus
all submultigraphs of $G$ that are isomorphic to a Shannon multigraph
of degree~$d$ are vertex-disjoint.

Let $G'$ be the multigraph that arises from $G$ by deleting every
edge from $G$ that belongs to a submultigraph of $G$ that is
isomorphic to a Shannon multigraph of degree~$d$. Then, $\Delta(G')
\le \Delta(G) = d$ and, by construction, $G'$ does not contain a
submultigraph of $G$ that is isomorphic to a Shannon multigraph of
degree~$d$. Since $xv \in E(G')$, the multigraph $G'$ has at least
one edge. By Shannon's Theorem and Vizing's Theorem, we deduce that
$\chi'(G) < \left\lfloor 3d/2 \right\rfloor$.

Let $\cC'$ be a $\chi'(G')$-edge coloring of the edges of $G'$. By
construction, every submultigraph of the connected multigraph $G$
that is isomorphic to a Shannon multigraph of degree~$d$ contains
exactly one vertex that is incident with an edge of $G'$. Since
$\chi'(G) < \left\lfloor 3d/2 \right\rfloor$, the coloring $\cC'$ can
therefore be extended to a $\left\lfloor 3d/2 \right\rfloor$-edge
coloring $\cC^*$ of $G$. Since $\cC'$ colors the edges of $G'$ with
fewer than $\left\lfloor 3d/2 \right\rfloor$ colors, $\cC^*$ is a
$\chi'(G)$-edge coloring of the edges of $G$ with at least two edge
color classes having different cardinality. This contradicts our
earlier observation that the edge color classes in every
$\chi'(G)$-edge coloring have the same cardinality. Therefore, $M$ is
not a proper submultigraph of the connected multigraph $G$, implying
that $G = M$. Hence if $d$ is odd, then $G$ is a Shannon multigraph
of degree $d$.

Conversely, if $G$ is a Shannon multigraph of degree $d$, then $m =
\left\lfloor 3d/2 \right\rfloor$
and $\alpha'(G) = 1$, implying that $\alpha'(G) = m/\left\lfloor 3d/2
\right\rfloor$.~\qed

\medskip
We close this section by recalling Hall's Matching Theorem  due to
K\"{o}nig~\cite{Ko31} and Hall~\cite{Ha35}.

\begin{unnumbered}{Hall's Matching Theorem.}
Let $G$ be a bipartite graph with partite sets $X$ and $Y$. Then $X$
can be matched to a subset of $Y$ if and only if $|N(S)| \ge |S|$ for
every nonempty subset $S$ of $X$.
\end{unnumbered}

\section{Proof of Main Result}
\label{S:proof}

We shall need the following properties of special hypergraphs defined in Section~\ref{S:special}.

\begin{ob}
Let $k \ge 2$ and let $H = E_k$ or $H = T_k$ and let $H$ have $n$
vertices and $m$ edges. Then the following holds. \\
\indent {\rm (a)} If $H = E_k$, then $\tau(H) = 1$. \\
\indent {\rm (b)} If $H = T_k$, then $\tau(H) = 2$. \\
\indent {\rm (c)} $\tau(H) = ( n + \left\lfloor \frac k2\right\rfloor
m ) / \left\lfloor \frac{3k}2 \right\rfloor$. \\
\indent {\rm (d)} Every vertex in $H$ belongs to some $\tau(H)$-set.
\label{snail}
\end{ob}

We are now in a position to prove our main result. Recall the
statement of Theorem~\ref{t:main}.

\noindent \textbf{Theorem~\ref{t:main}}. \emph{For $k = 2$
or $k \ge 4$, let $H$ be a connected
$k$-uniform hypergraph on $n$ vertices and $m$ edges. Then,
\[
\tau(H) \le \frac{n + \left\lfloor \frac k2\right\rfloor m}{\left\lfloor \frac{3k}2\right\rfloor}
\]
\\
\noindent with equality if and only if  $H = E_k$ or $H = T_k$. }

\noindent \textbf{Proof.} The upper bound on $\tau(H)$ is a
restatement of the Chv\'{a}tal-McDiarmid Theorem. We only need prove
that $\tau(H) = ( n + \left\lfloor \frac k2\right\rfloor m ) /
\left\lfloor \frac{3k}2 \right\rfloor$ if and only if  $H = E_k$ or $H = T_k$. If $H = E_k$ or $H = T_k$, then by
Observation~\ref{snail}, $\tau(H) = ( n + \left\lfloor \frac
k2\right\rfloor m ) / \left\lfloor \frac{3k}2 \right\rfloor$, as
desired.

To prove the converse, suppose that $\tau(H) = ( n + \left\lfloor
\frac k2\right\rfloor m ) / \left\lfloor \frac{3k}2 \right\rfloor$,
where $k = 2$ or $k \ge 4$. We proceed by induction on the order~$n$
to show that $H = E_k$ or $H = T_k$. If $m = 0$, then $\tau(H) = 0
< n / \left\lfloor \frac{3k}2 \right\rfloor$, a contradiction. Hence
$m \ge 1$, and so $n \ge k$. If $n = k$, then $H = E_k$,
and we are done. This establishes the base case. Let $n \ge k+1$ and
let $H$ be a connected $k$-uniform hypergraph on $n$ vertices and $m$
edges, and assume that the desired result holds for all connected
$k$-uniform hypergraph on fewer than $n$ vertices.

\begin{unnumbered}{Claim~A}
$\delta(H) \ge 1$.
\end{unnumbered}
\textbf{Proof.} Suppose that $\delta(H) = 0$. Let $F$ be obtained
from $H$ by deleting all isolated vertices. Let $F$ have $\nF$
vertices and $\mF$ edges. Then, $\nF \le n - 1$ and $\mF = m$. Every
transversal in $H'$ is a transversal in $H$, and so $\tau(H) \le
\tau(H')$. By the Chv\'{a}tal-McDiarmid Theorem, we have that
\[
\tau(H) \le \tau(H') \le
\frac{\nF + \left\lfloor \frac k2\right\rfloor \mF}{\left\lfloor \frac{3k}2\right\rfloor}
 <
\frac{n + \left\lfloor \frac k2\right\rfloor m}{\left\lfloor
\frac{3k}2\right\rfloor},
\]
a contradiction. Hence, $\delta(H) \ge 1$.~\smallqed

\medskip
Let $v$ be a vertex of maximum degree~$\Delta(H)$ in $H$ and let
$H'=H-v$ have $n'$ vertices and $m'$ edges. Then, $H'$ is a
$k$-uniform hypergraph. Every transversal in $H'$ can be extended to
a transversal in $H$ by adding to it the vertex $v$, and so $\tau(H)
\le \tau(H') + 1$. Recall that $k = 2$ or $k \ge 4$.

\begin{unnumbered}{Claim~B}
If $k$ is even, then $\Delta(H)\le 2$, while if $k$ is odd, then
$\Delta(H) \le 3$.
\end{unnumbered}
\textbf{Proof.} Suppose first that $k$ is even and $\Delta(H) \ge 3$.
Then, $n'\le n - 1$ and $m' \le m - 3$. Since $k$ is even, we have by
the Chv\'{a}tal-McDiarmid Theorem that
\[
\tau(H) \le \tau(H') + 1 \le
\frac{n' + \left\lfloor \frac k2\right\rfloor m'}{\left\lfloor \frac{3k}2\right\rfloor} + 1
\le \frac{n + \left\lfloor \frac k2\right\rfloor m - 1}{\left\lfloor
\frac{3k}2\right\rfloor}
< \frac{n + \left\lfloor \frac k2\right\rfloor m}{\left\lfloor
\frac{3k}2\right\rfloor},
\]
a contradiction. Hence if $k$ is even, then $\Delta(H) \le 2$.
Suppose next that $k$ is odd and $\Delta(H) \ge 4$. Then, $n' \le n -
1$ and $m' \le m - 4$. Since $k \ge 5$ is odd, we have by the
Chv\'{a}tal-McDiarmid Theorem that
\[
\tau(H) \le \tau(H') + 1 \le
\frac{n' + \left\lfloor \frac k2\right\rfloor m'}{\left\lfloor \frac{3k}2\right\rfloor} + 1
\le \frac{(n-1) +  \left\lfloor \frac k2\right\rfloor (m-4)}{\left\lfloor \frac{3k}2\right\rfloor}+1
< \frac{n + \left\lfloor \frac k2\right\rfloor m}{\left\lfloor
\frac{3k}2\right\rfloor},
\]

\noindent a contradiction. Hence if $k$ is odd, then $\Delta(H) \le
3$.~\smallqed

\begin{unnumbered}{Claim~C}
If $k = 2$, then $H$ is a generalized triangle $T_2$.
\end{unnumbered}
\textbf{Proof.} Suppose that $k=2$, and so $H$ is a graph and
$\tau(H) = (n+m)/3$. By Claim~A and Claim~B, we have that $\delta(H)
\ge 1$ and $\Delta(H) \le 2$. Thus, $H$ is a path or a cycle. If $H$
is a path on $n \ge 2$ vertices, then $(2n-1)/3 = (n+m)/3 = \tau(H) =
\lfloor n/2 \rfloor$, implying that $n = 2$ and $H = E_2$.
However this contradicts the fact that $n \ge k+1$. Hence, $H$ is a
cycle on $n \ge 3$ vertices. Thus, $2n/3 = (n+m)/3 = \tau(H) = \lceil
n/2 \rceil$, implying that $n = 3$ and $H$ is a generalized triangle $T_2$.~\qed

\medskip
In what follows we may assume that $k \ge 4$, for otherwise the
desired result follows by Claim~C.

\begin{unnumbered}{Claim~D}
If $\Delta(H) \le 2$, then $H = T_k$.
\end{unnumbered}
\textbf{Proof.} Suppose that $\Delta(H) \le 2$. For $i=1,2$, let
$n_i$ be the number of vertices of degree~$i$ in $H$. By Claim~A,
$\delta(H) \ge 1$ and so $n_1 + n_2 = n$. By the $k$-uniformity of
$H$ we have that $n_1 + 2n_2=km$, or, equivalently, $n_2 =km-n$. We
now consider the multigraph $G$ whose vertices are the edges of $H$
and whose edges correspond to the $n_2$ vertices of degree~$2$ in
$H$: if a vertex of $H$ is contained in the edges $e$ and $f$ of $H$,
then the corresponding edge of $G$ joins vertices $e$ and $f$ of $G$.
Since $H$ is $k$-uniform and $\Delta(H) \le 2$, the maximum degree in
$G$ is at most~$k$. Further since $H$ is connected, so too is $G$.

Let $M$ be a maximum matching in $G$, and so by
Theorem~\ref{matching}, $|M| = \alpha'(G) \ge n_2/\left\lfloor 3k/2
\right\rfloor$. Let $S$ be the set of vertices of $H$ that correspond
to the set of edges $M$ in $G$. Then, $S$ is an independent set in
$H$ and every vertex in $S$ has degree~$2$ in $H$. By the maximality
of $M$, we note that the set of edges in $H$ that do not intersect
$S$ are vertex-disjoint. Let $S'$ be a set of vertices in $H$ that
consists of exactly one vertex from every edge of $H$ that does not
intersect $S$. Then, $|S'| = m - 2|S|$ and the set $S \cup S'$ is a
transversal in $H$. Thus, $\tau(H) \le |S| + |S'| = m - |S| = m -
|M|$. Hence,
\[
\frac{n + \left\lfloor \frac k2\right\rfloor m}{\left\lfloor
\frac{3k}2\right\rfloor}
= \tau(H) \le m-|M|\le m-\frac{n_2}{\left\lfloor\frac{3k}2 \right\rfloor}
=m-\frac{km-n}{\left\lfloor \frac{3k}2\right\rfloor}
=\frac{\left\lfloor \frac k2\right\rfloor m+n}{\left\lfloor \frac{3k}2\right\rfloor}.
\]

\noindent Consequently, we must have equality throughout the above
inequality chain. In particular, $\alpha'(G) = |M| = n_2/\left\lfloor
3k/2 \right\rfloor$. Thus by Theorem~\ref{matching}, either $n_2 = 0$
or $G$ is a Shannon multigraph of degree~$k$. If $n_2 = 0$, then $n =
n_1$, implying by the connectivity of $H$ that $H = E_k$.
However this contradicts the fact that $n \ge k+1$. Hence, $G$ is a
Shannon multigraph of degree~$k$, implying that $H$ is a generalized triangle $T_k$.~\smallqed

\medskip
By Claim~D, if $\Delta(H) \le 2$, then $H$ is a generalized triangle $T_k$,
and we are done. Hence we may assume in what follows that $\Delta(H)
\ge 3$. By Claim~B, $k \ge 5$ is odd and $\Delta(H) = 3$.
We now prove a series of claims that culminate in a
contradiction.\footnote{We remark that if we allow $k=3$, then it is
indeed possible that $\Delta(H) = 3$. The current proof technique
therefore fails in this special case when $k = 3$ since we are then
unable to associate a multigraph with the hypergraph $H$ as is done
in the proof of Claim~D. However as remarked earlier, the special
case when $k=3$ has fortunately been handled in~\cite{HeYe08}.}

\begin{unnumbered}{Claim~E}
The following hold in the hypergraph $H$. \\
\indent {\rm (a)} $\tau(H)=\tau(H')+1$.
\\
\indent {\rm (b)} $n' = n-1$.
\\
\indent {\rm (c)} Every component of $H'$ is either $E_k$ or $T_k$.
\end{unnumbered}
\textbf{Proof.} Since $\Delta(H) = 3$, we note that $n' \le n - 1$
and $m' = m - 3$. Since $k$ is odd, we have by the
Chv\'{a}tal-McDiarmid Theorem that
\[
\tau(H) \le \tau(H') + 1 \le
\frac{n' + \left\lfloor \frac k2\right\rfloor m'}{\left\lfloor \frac{3k}2\right\rfloor} + 1
\le \frac{(n-1) +  \left\lfloor \frac k2\right\rfloor (m-3)}{\left\lfloor \frac{3k}2\right\rfloor}+1
= \frac{n + \left\lfloor \frac k2\right\rfloor m}{\left\lfloor
\frac{3k}2\right\rfloor},
\]

\noindent Since $\tau(H) = ( n + \left\lfloor \frac k2\right\rfloor m
) / \left\lfloor \frac{3k}2 \right\rfloor$, we must have equality
throughout the above inequality chain, implying that
$\tau(H)=\tau(H')+1$, $\tau(H') = ( n' + \left\lfloor \frac
k2\right\rfloor m' ) / \left\lfloor \frac{3k}2 \right\rfloor$ and
$n'=n-1$. Applying the inductive hypothesis to every component of
$H'$, we have that every component of $H'$ is either $E_k$ or $T_k$.~\smallqed

\medskip
By Claim~E(c) every component of $H'$ is either $E_k$ or $T_k$.
By Observation~\ref{snail}, every component of $H'$
that is $E_k$ or $T_k$ contributes $1$ or $2$, respectively, to $\tau(H')$.

Let $e_1, e_2, e_3$ be the three edges that contain the vertex~$v$ in
$H$ and let $E_v = \{e_1,e_2,e_3\}$. By Claim~E(b), $n' = n-1$, which
implies that $|V(e) \cap V(H')| = k-1$ for each edge $e \in E_v$.

\begin{unnumbered}{Claim~F}
Let $C$ be a component of $H'$ that is a generalized triangle $T_k$. If $|V(C) \cap V(e_1)|\ge 2$,
$|V(C)\cap V(e_2)| \ge 2$ and $|V(C)\cap V(e_3)| \le k-2$, then
$|V(C)\cap V(e_1)| + |V(C)\cap V(e_2)| \le (k+1)/2$.
\end{unnumbered}
\textbf{Proof.} Assume, to the contrary, that $|V(C) \cap V(e_1)| +
|V(C)\cap V(e_2)| > (k+1)/2$. Since $|V(C) \cap V(e_3)| \le k-2$,
there is a vertex $u_3 \in V(e_3) \setminus (V(C)\cup\{v\})$. If
there is a vertex $u_1 \in V(C) \cap V(e_1) \cap V(e_2)$, then by
Observation~\ref{snail}(d) there is a $\tau(H')$-set $T$ that
contains both $u_1$ and $u_3$. Since $\{u_1,u_3\}$ intersects all
three edges that contain~$v$ in $H$, the set $T$ is a transversal of
$H$, and so $\tau(H) \le |T| = \tau(H')$, contradicting Claim~E(a).
Hence, $V(C) \cap V(e_1) \cap V(e_2) = \emptyset$.

Since $|V(C) \cap V(e_1)|\ge 2$ and since there is a unique vertex of $C$ of degree~$1$ in $H'$, there is a vertex $u_1
\in V(C) \cap V(e_1)$ of degree~$2$ in $H'$. Let $f$ be the unique
edge of $C$ that does not contain~$u_1$. If there is a vertex $u_2
\in V(f) \cap V(e_2)$, then $\{u_1,u_2\}$ is a $\tau(C)$-set and, by
Observation~\ref{snail}(d), there is a $\tau(H')$-set $T$ that
contains the set $\{u_1,u_2,u_3\}$. Since $\{u_1,u_2,u_3\}$
intersects all three edges that contain~$v$ in $H$, the set $T$ is a
transversal of $H$, and so $\tau(H) \le |T| = \tau(H')$,
contradicting Claim~E(a). Hence, $V(f) \cap V(e_2) = \emptyset$.

Since $|V(C) \cap V(e_2)|\ge 2$ and since there is a unique vertex of $C$ of degree~$1$ in $H'$, there is a vertex $x_2
\in V(C) \cap V(e_2)$ of degree~$2$ in $H'$. By assumption,
$|V(C)\cap V(e_1)| + |V(C) \cap V(e_2)| > (k+1)/2$. As observed
earlier, the edges $e_1$ and $e_2$ do not intersect in $C$ and $V(f)
\cap V(e_2) = \emptyset$. Since $|V(C) \setminus V(f)| = (k+1)/2$,
there is therefore a vertex $x_1 \in V(f) \cap V(e_1)$. Thus the set
$\{x_1,x_2\}$ is a $\tau(C)$-set and, by Observation~\ref{snail}(d),
there is a $\tau(H')$-set $T$ that contains the set
$\{x_1,x_2,u_3\}$. Since $\{x_1,x_2,u_3\}$ intersects all three edges
that contain~$v$ in $H$, the set $T$ is a transversal of $H$, and so
$\tau(H) \le |T| = \tau(H')$, once again contradicting Claim~E(a).
Therefore, $|V(C)\cap V(e_1)| + |V(C)\cap V(e_2)| \le
(k+1)/2$.~\smallqed

\begin{unnumbered}{Claim~G}
$H'$ is disconnected.
\end{unnumbered}
\textbf{Proof.} Assume, to the contrary, that $H'$ is connected.
Since $\Delta(H)=3$, the $k$-uniformity of $H$ implies that $km =
\sum_{v \in V(H)} d(v) \le 3n$. By Claim~E(c), $H$ is either
$E_k$ or $T_k$. Suppose first that $H = E_k$. Then,
$n=k+1$ and $m=4$. However $k \ge 5$, and so $km = 4k\ge 3k+5
> 3k+3=3n$, a contradiction. Hence, $H = T_k$. Thus, $n =
3(k+1)/2$ and $m = 6$. However $k \ge 5$, and so $km = 6k = 9k/2 +
3k/2 \ge 9k/2 + 15/2 > 9k/2 + 9/2 = 3n$, once again producing a
contradiction. Therefore, $H'$ is disconnected.~\smallqed

\begin{unnumbered}{Claim~H}
$H'$ has at least three components.
\end{unnumbered}
\textbf{Proof.} Assume, to the contrary, that $H'$ has at most two
components. Then by Claim~G, the hypergraph $H'$ has exactly two
components which we call $C_1$ and $C_2$. As observed earlier, $|V(e)
\cap V(H')| = k-1$ for each edge $e \in E_v$. Renaming the components
$C_1$ and $C_2$ if necessary, we may assume that
\begin{equation}
\sum_{e \in E_v} |V(C_1)\cap V(e)| \ge \frac{3}{2}(k-1) \ge \sum_{e \in E_v} |V(C_2) \cap V(e)|
\label{Eq1}
\end{equation}
and that if we have equality throughout the Inequality
Chain~(\ref{Eq1}), then $V(C_2)$ intersects as least as many edges of
$E_v$ as $V(C_1)$ does. Since $H$ is connected, the vertex $v$ is
adjacent in $H$ to a vertex from $V(C_1)$ and to a vertex from
$V(C_2)$.

\begin{unnumbered}{Claim~H.1}
$C_1 = T_k$.
\end{unnumbered}
\textbf{Proof.} Assume, to the contrary, that $C_1 = E_k$, and so $C_1$ has $k$ vertices. By our choice of $C_1$, $\sum_{e \in E_v}
|V(C_1)\cap V(e)| \ge 3(k-1)/2$. Since $k \ge 5$, we have that
$3(k-1)/2 > k$. Hence by the pigeonhole principle, at least one
vertex, $u_1$ say, of $C_1$ is contained in two edges of $E_v$, and
so $u_1$ has degree~$3$ in $H$. Renaming the edges in $E_v$ if
necessary, we may assume that $u_1 \in V(e_1) \cap V(e_2)$.

If the edge $e_3$ intersects $V(C_2)$, then let $u_3 \in V(C_2) \cap
V(e_3)$. By Observation~\ref{snail}(d), there is a $\tau(H')$-set $T$
that contains the set $\{u_1,u_3\}$. Since $\{u_1,u_3\}$ intersects
all three edges that contain~$v$ in $H$, the set $T$ is a transversal
of $H$, and so $\tau(H) \le |T| = \tau(H')$, contradicting
Claim~E(a). Hence the edge $e_3$ does not intersect $V(C_2)$. Thus,
$V(e_3) \setminus \{v\} \subset V(C_1)$.

Suppose that both edges $e_1$ and $e_2$ intersect $V(C_2)$. If all
vertices of $V(C_1) \cap V(e_3)$ have degree~$2$ in $H$, then $V(e_3)
\setminus \{v\} = V(C_1) \setminus \{u_1\}$ and $V(C_1) \cap V(e_1) =
\{u_1\} = V(C_1) \cap V(e_2)$. Thus,
$3(k-1)/2 \le \sum_{e \in E_v} |V(C_1)\cap V(e)| = k+1$, and so $k
\le 5$. Consequently, $k = 5$ and we have equality throughout the
Inequality Chain~(\ref{Eq1}). But then all three edges in $E_v$
intersect $V(C_1)$ but only two edges in $E_v$ intersect $V(C_2)$,
contradicting our choice of $C_1$ and $C_2$. Therefore there is a
vertex $x_1 \in V(C_1) \cap V(e_3)$ that has degree~$3$ in $H$.
Renaming the edges $e_1$ and $e_2$, if necessary, we may assume that
$x_1 \in V(e_1)$. By assumption, the edge $e_2$ intersect $V(C_2)$.
Let $x_2 \in V(C_2) \cap V(e_2)$. By Observation~\ref{snail}(d),
there is a $\tau(H')$-set $T$ that contains the set $\{x_1,x_2\}$.
Since $\{x_1,x_2\}$ intersects all three edges that contain~$v$ in
$H$, the set $T$ is a transversal of $H$, and so $\tau(H) \le |T| =
\tau(H')$, contradicting Claim~E(a). Hence, at most one of $e_1$ and $e_2$ intersects
$V(C_2)$.

Hence renaming $e_1$ and $e_2$, if necessary, we may assume that
$V(e_1) \setminus \{v\} \subset V(C_1)$. By the pigeonhole principle,
there is a vertex $w_1 \in V(C_1) \cap V(e_1) \cap V(e_3)$. Since $H$
is connected, the edge $e_2$ intersects $V(C_2)$. Let $w_2 \in V(C_2)
\cap V(e_2)$. By Observation~\ref{snail}(d), there is a
$\tau(H')$-set $T$ that contains the set $\{w_1,w_2\}$. Since
$\{w_1,w_2\}$ intersects all three edges that contain~$v$ in $H$, the
set $T$ is a transversal of $H$, and so $\tau(H) \le |T| = \tau(H')$,
contradicting Claim~E(a). Therefore, $C_1$ is a generalized triangle
$T_k$.~\smallqed

\medskip
By Claim~H.1, the component $C_1$ is a generalized triangle
$T_k$. Renaming the
edges $e_1, e_2, e_3$ if necessary, we may assume that
\[
|V(C_1)\cap V(e_1)| \ge |V(C_1)\cap V(e_2)| \ge |V(C_1)\cap
V(e_3)|,
\]
which implies that
\[
|V(C_1)\cap V(e_3)| \le \frac13 \sum_{e \in E_v} |V(C_1)\cap V(e)|.
\]
Therefore,
\[
|V(C_1)\cap V(e_1)|+|V(C_1)\cap V(e_2)| \ge \frac23 \sum_{e \in E_v} |V(C_1)\cap V(e)| \ge k-1 > \frac{1}{2}(k+1).
\]

If $e_3$ does not intersect $V(C_2)$, then neither do the edges $e_1$
and $e_2$, implying that $H$ is disconnected, a contradiction. Hence,
$e_3$ intersect $V(C_2)$, and so $|V(C_1) \cap V(e_3)| \le k-2$. If
$|V(C_1) \cap V(e_2)| \ge 2$, then $|V(C_1) \cap V(e_1)| \ge 2$. But
then we contradict Claim~F. Therefore, $|V(C_1)\cap V(e_2)| \le 1$,
and so $|V(C_1)\cap V(e_3)| \le 1$. Now by our choice of $C_1$,
\[
k+1 \le \frac{3}{2}(k-1) \le \sum_{e \in E_v} |V(C_1)\cap V(e)| \le (k-1)+1+1=k+1.
\]
\indent Consequently, we must have equality throughout the above
inequality chain. In particular, $\sum_{e \in E_v} |V(C_1)\cap V(e)|
= 3(k-1)/2$, $|V(C_1) \cap V(e_1)|=k-1$ and $|V(C_1) \cap V(e_2)| =
|V(C_1)\cap V(e_3)| = 1$. But then we have equality throughout the Inequality
Chain~(\ref{Eq1}) and all three edges in $E_v$ intersect
$V(C_1)$ but only two edges in $E_v$ intersect $V(C_2)$,
contradicting our choice of $C_1$ and $C_2$. Therefore, $H'$ has at
least three components. This completes the proof of
Claim~H.~\smallqed

\medskip
We now return to the proof of Theorem~\ref{t:main}. By Claim~H, the
hypergraph $H'$ has at least three components. Let $F$ be a bipartite
graph with partite sets $V_1$ and $V_2$, where $V_1 = E_v =
\{e_1,e_2,e_3\}$ and where the vertices in $V_2$ correspond to the
components of $H'$. Further the edge set of $F$ is defined as
follows: If an edge $e \in E_v$ intersects a component $C$ of $H'$ in
$H$, then the vertex $e \in V_1$ is adjacent to the vertex $C \in
V_2$ in $F$.

Since $H'$ has at least three components, $|V_2| \ge 3$. Since $H$ is
connected, every component in $H'$ has a nonempty intersection with
at least one edge in $E_v$, and so every vertex in $V_2$ has degree
at least~$1$ in $F$ and $N_F(V_1) = V_2$. Thus if $S = V_1$, then
$|N_F(S)| = |V_2| \ge 3 = |S|$. Since every edge $e \in E_v$
intersects at least one component of $H'$ in $H$, every vertex in
$V_1$ has degree at least~$1$ in $F$. Thus if $S \subset V_1$ and
$|S| = 1$, then $|N_F(S)| \ge |S|$. Hence by Hall's Matching Theorem,
either $V_1$ can be matched to a subset of $V_2$ in $F$ or $|N_F(S)|
< |S|$ for some subset $S \subset V_1$ with $|S| = 2$.

Suppose that $V_1$ can be matched to a subset of $V_2$ in $F$. Let
$M_F$ be such a matching in $F$. We now name the components in $H'$
so that $M_F = \{e_1C_1,e_2C_2,e_3C_3\}$. Hence for $i \in
\{1,2,3\}$, the edge $e_i$ intersects the component $C_i$ of $H'$ in
$H$. For $i \in \{1,2,3\}$, let $u_i \in V(C_i) \cap V(e_i)$. By
Observation~\ref{snail}(d), there is a $\tau(H')$-set $T$ that
contains the set $\{u_1,u_2,u_3\}$. Since $\{u_1,u_2,u_3\}$
intersects all three edges that contain~$v$ in $H$, the set $T$ is a
transversal of $H$, and so $\tau(H) \le |T| = \tau(H')$,
contradicting Claim~E(a). Therefore, $|N_F(S)| < |S|$ for some subset
$S \subseteq V_1$ with $|S| = 2$.

Renaming the edges in $E_v$ if necessary, we may assume that $S =
\{e_1,e_2\}$. Thus in $H$ we have that $V(e_1), V(e_2) \subseteq V(C)
\cup\{v\}$ for some component $C$ of $H'$. Since $H$ is connected,
the edge $e_3$ intersects every component of $H'$ different from $C$
in $H$. Thus, $|V(C) \cap V(e_1)| = k-1$, $|V(C) \cap V(e_2)| = k -
1$ and $|V(C) \cap V(e_3)| \le k-3$. If $C$ is a generalized triangle
$T_k$, then we contradict Claim~F. Hence, $C = E_k$.

Let $C'$ be an arbitrary component of $H'$ different from $C$, and
let $u_3 \in V(C') \cap V(e_3)$. Since $\sum_{i=1}^{2} |V(C) \cap
V(e_i)| = 2(k-1) > k$, by the pigeonhole principle at least one vertex,
$u_1$ say, of $C$ is contained in both edges $e_1$ and $e_2$. By
Observation~\ref{snail}(d), there is a $\tau(H')$-set $T$ that
contains the set $\{u_1,u_3\}$. Since $\{u_1,u_3\}$ intersects all
three edges that contain~$v$ in $H$, the set $T$ is a transversal of
$H$, and so $\tau(H) \le |T| = \tau(H')$, contradicting Claim~E(a).
This completes the proof of Theorem~\ref{t:main}.~\qed

\section{Acknowledgements}

Research of the first author is supported in part by the South
African National Research Foundation and the University of
Johannesburg, and research of the second author is supported by the
Deutsche Forschungsgemeinschaft (GZ: LO 1758/1-1).

\end{document}